\documentclass[reqno,a4paper]{amsart}
\usepackage{amssymb,setspace}
\usepackage{ifpdf}
\ifpdf
  \usepackage[hyperindex]{hyperref}
\else
  \expandafter\ifx\csname dvipdfm\endcsname\relax
    \usepackage[hypertex,hyperindex]{hyperref}
  \else
    \usepackage[dvipdfm,hyperindex]{hyperref}
  \fi
\fi
\allowdisplaybreaks[4]
\theoremstyle{plain}
\newtheorem{thm}{Theorem}
\theoremstyle{remark}
\newtheorem{rem}{Remark}
\DeclareMathOperator{\td}{d\mspace{-2mu}}

\begin{document}

\title[Monotonic functions involving polygamma functions]
{Complete monotonicity of a family of functions involving the tri- and tetra-gamma functions}

\author[F. Qi]{Feng Qi}
\address{Department of Mathematics, College of Science, Tianjin Polytechnic University, Tianjin City, 300160, China}
\email{\href{mailto: F. Qi <qifeng618@gmail.com>}{qifeng618@gmail.com}, \href{mailto: F. Qi <qifeng618@hotmail.com>}{qifeng618@hotmail.com}, \href{mailto: F. Qi <qifeng618@qq.com>}{qifeng618@qq.com}}
\urladdr{\url{http://qifeng618.wordpress.com}}

\begin{abstract}
The psi function $\psi(x)$ is defined by $\psi(x)=\frac{\Gamma'(x)}{\Gamma(x)}$ and $\psi^{(i)}(x)$ for $i\in\mathbb{N}$ denote polygamma functions, where $\Gamma(x)$ is the gamma function. In this paper, we prove that the function
$$
[\psi'(x)]^2+\psi''(x)-\frac{x^2+\lambda x+12}{12x^4(x+1)^2}
$$
is completely monotonic on $(0,\infty)$ if and only if $\lambda\le0$, and so is its negative if and only if $\lambda\ge4$. From this, some inequalities are refined and sharpened.
\end{abstract}



\thanks{This paper was typeset using \AmS-\LaTeX}

\maketitle

\section{Introduction}

We recall from~\cite[Chapter~XIII]{mpf-1993} and~\cite[Chapter~IV]{widder} that a function
$f$ is said to be completely monotonic on an interval $I$ if $f$ has derivatives of all
orders on $I$ and
\begin{equation}\label{CM-dfn}
0\le(-1)^{n}f^{(n)}(x)<\infty
\end{equation}
for $x\in I$ and $n\ge0$. The class of completely monotonic functions may be characterized by the famous Bernstein\nobreakdash-Widder Theorem~\cite[p.~161, Theorem~12b]{widder} which reads that a necessary and sufficient condition that $f(x)$ should be completely monotonic for $0<x<\infty$ is that
\begin{equation} \label{berstein-1}
f(x)=\int_0^\infty e^{-xt}\td\alpha(t),
\end{equation}
where $\alpha(t)$ is non-decreasing and the integral converges for $0<x<\infty$.
\par
We also recall that the classical Euler's gamma function $\Gamma(x)$ is defined by
\begin{equation}\label{gamma-dfn}
\Gamma(x)=\int^\infty_0t^{x-1} e^{-t}\td t,\quad x>0.
\end{equation}
The logarithmic derivative of $\Gamma(x)$, denoted by $\psi(x)=\frac{\Gamma'(x)}{\Gamma(x)}$, is called the psi or di-gamma function, and the derivatives $\psi^{(i)}(x)$ for $i\in\mathbb{N}$ are respectively called the polygamma functions. In particular, the functions $\psi'(x)$ and $\psi''(x)$ are called the tri- and tetra-gamma functions.
\par
In~\cite[p.~208, (4.39)]{forum-alzer}, it was established that the one-sided inequality
\begin{equation}\label{psi'psi''-lower}
[\psi'(x)]^2+\psi''(x)>\frac{p(x)}{900x^4(x+1)^{10}}
\end{equation}
holds for $x>0$, where
\begin{equation}
\begin{split}\label{p(x)-dfn}
p(x)&=75x^{10}+900x^9+4840x^8+15370x^7+31865x^6+45050x^5\\
&\quad+44101x^4+29700x^3+13290x^2+3600x+450.
\end{split}
\end{equation}
By the same technique as in ~\cite[p.~208, (4.39)]{forum-alzer}, we can easily obtain
\begin{equation}\label{psi'psi''-upper-weaker}
  [\psi'(x)]^2+\psi''(x)<\frac{36+180 x+408 x^2+504 x^3+352 x^4+132 x^5+21 x^6}{36 x^4 (1+x)^6}
\end{equation}
for $x>0$.
\par
In \cite{x-4-di-tri-gamma-p(x)-Slovaca.tex}, the function
\begin{equation}\label{psi'psi''-function}
[\psi'(x)]^2+\psi''(x)-\frac{p(x)}{900x^4(x+1)^{10}}
\end{equation}
was proved to be completely monotonic on $(0,\infty)$.
\par
In \cite{x-4-di-tri-gamma-upper-lower-combined.tex}, the following results were obtained:
\begin{enumerate}
\item
The two-sided inequality
\begin{equation}\label{psi'psi''-upper}
\frac{x^2+12}{12x^4(x+1)^2}<[\psi'(x)]^2+\psi''(x)<\frac{x+12}{12x^4(x+1)}
\end{equation}
holds on $(0,\infty)$.
\item
The functions
\begin{equation}\label{x-4-upper-f(x)}
f(x)=[\psi'(x)]^2+\psi''(x)-\frac{x^2+12}{12x^4(x+1)^2}
\end{equation}
and
\begin{equation}\label{x-4-upper-g(x)}
g(x)=\frac{x+12}{12x^4(x+1)}-\bigl\{[\psi'(x)]^2+\psi''(x)\bigr\}
\end{equation}
are completely monotonic on $(0,\infty)$.
\end{enumerate}
\par
The inequalities \eqref{psi'psi''-lower} and~\eqref{psi'psi''-upper-weaker} and the ones in~\eqref{psi'psi''-upper} are not included each other.
\par
For more information on results related to the function $[\psi'(x)]^2+\psi''(x)$, please refer to~\cite{notes-best-simple-open-jkms.tex, psi-proper-fraction-degree-two.tex, Yang-Fan-2008-Dec-simp.tex, Guo-Qi-Srivasta-Unique.tex, BAustMS-5984-RV.tex, notes-best-simple-cpaa.tex, AAM-Qi-09-PolyGamma.tex, SCM-2012-0142.tex, x-4-di-tri-gamma-p(x)-Slovaca.tex}, the expository and survey article~\cite{bounds-two-gammas.tex, Wendel-Gautschi-type-ineq-Banach.tex} and the literature listed therein.
\par
The aim of this paper is to establish necessary and sufficient conditions on $\lambda\in\mathbb{R}$ for the function
\begin{equation}\label{tetra-square+poly-2}
f_\lambda(x)=[\psi'(x)]^2+\psi''(x)-\frac{x^2+\lambda x+12}{12x^4(x+1)^2}
\end{equation}
to be completely monotonic on $(0,\infty)$.
\par
Our main results may be stated as the following theorem.

\begin{thm}\label{x-4-upper-improve-thm}
Let $\lambda\in\mathbb{R}$.
\begin{enumerate}
\item
The function $f_\lambda(x)$ defined by~\eqref{tetra-square+poly-2} is completely monotonic on $(0,\infty)$ if and only if $\lambda\le0$;
\item
The function $-f_\lambda(x)$ is completely monotonic on $(0,\infty)$ if and only if $\lambda\ge4$;
\item
The double inequality
\begin{equation}\label{x-4-upper-improve-ineq}
\frac{x^2+\mu x+12}{12x^4(x+1)^2}<[\psi'(x)]^2+\psi''(x)<\frac{x^2+\nu x+12}{12x^4(x+1)^2}
\end{equation}
holds on $(0,\infty)$ if and only if $\mu\le0$ and $\nu\ge4$.
\end{enumerate}
\end{thm}

In next section we supply several proofs for Theorem~\ref{x-4-upper-improve-thm}. In the final section we derive some corollaries and pose a double inequality of $[\psi'(x)]^2+\psi''(x)$ on $(0,\infty)$.

\section{Proofs of Theorem~\ref{x-4-upper-improve-thm}}

In this section we provide several proofs for Theorem~\ref{x-4-upper-improve-thm} by different approaches.

\begin{proof}[First proof of Theorem~\ref{x-4-upper-improve-thm}]
By the recursion formula
\begin{equation}\label{psisymp4}
\psi^{(n-1)}(x+1)=\psi^{(n-1)}(x)+\frac{(-1)^{n-1}(n-1)!}{x^n}
\end{equation}
for $x>0$ and $n\in\mathbb{N}$, see \cite[pp.~258 and 260, 6.3.5 and 6.4.6]{abram}, we have
\begin{multline}
f_\lambda(x)-f_\lambda(x+1)=\bigl[\psi'(x)-\psi'(x+1)\bigr]\bigl[\psi'(x)+\psi'(x+1)\bigr]\\
\begin{aligned}
&\quad+\bigl[\psi''(x)-\psi''(x+1)\bigr] -\biggl[\frac{x^2+\lambda x+12}{12x^4(x+1)^2} -\frac{(x+1)^2+\lambda(x+1)+12}{12(x+1)^4(x+2)^2}\biggr]\\
&=\frac1{x^2}\biggl[2\psi'(x)-\frac1{x^2}\biggr]-\frac{2}{x^3}
-\biggl[\frac{x^2+\lambda x+12}{12x^4(x+1)^2} -\frac{(x+1)^2+\lambda(x+1)+12}{12(x+1)^4(x+2)^2}\biggr]\\
&=\frac2{x^2}\biggl\{\psi'(x)-\frac1{2x^2}-\frac1{x} -\frac{x^2+\lambda x+12}{24x^2(x+1)^2} +\frac{x^2[(x+1)^2+\lambda(x+1)+12]}{24(x+1)^4(x+2)^2}\biggr\}\\
&=\frac2{x^2}\biggl[\psi'(x)-\frac{1}{x^2}-\frac{\lambda}{24 x}+\frac{37-2 \lambda}{6 (x+2)} +\frac{9\lambda-172}{24(x+1)}+\frac{28-\lambda}{8 (x+1)^2}\\
&\quad+\frac{13-\lambda}{6(x+2)^2}+\frac{\lambda-48}{24 (x+1)^3}+\frac{1}{2(x+1)^4}\biggr]\\
&\triangleq \frac2{x^2}h_\lambda(x).
\end{aligned}
\end{multline}
Using the formula
\begin{equation}\label{gam}
\frac1{x^r}=\frac1{\Gamma(r)}\int_0^\infty t^{r-1}e^{-xt}\td t
\end{equation}
for $r>0$ and $x>0$, see \cite[p.~255, 6.1.1]{abram}, and the integral representations
\begin{equation} \label{psin}
\psi ^{(n)}(x)=(-1)^{n+1}\int_{0}^{\infty}\frac{t^{n}}{1-e^{-t}}e^{-xt}\td t
\end{equation}
for $n\in\mathbb{N}$ and $x\in(0,\infty)$, see \cite[p.~260, 6.4.1]{abram}, yields
\begin{align*}
h_\lambda(x)&=\int_0^\infty\biggl(\frac{t}{1-e^{-t}}-t-\frac\lambda{24}+\frac{37-2\lambda}6e^{-2t} +\frac{13-\lambda}6te^{-2t} \\
&\quad+\frac{9\lambda-172}{24}e^{-t}+\frac{28-\lambda}8te^{-t}+\frac{\lambda-48}{48}t^2e^{-t} +\frac1{12}t^3e^{-t}\biggr)e^{-xt}\td t\\
&=\frac1{48}\int_0^\infty\biggl[\frac{4P(t)}{Q(t)}-\lambda\biggr]Q(t)e^{-(x+2)t}\td t,
\end{align*}
where
\begin{equation*}
P(t)=\frac{e^{2 t} \bigl(t^3-12 t^2+54 t-86\bigr)-e^t \bigl(t^3-12 t^2+16 t-160\bigr)-26 t-74}{e^t-1}
\end{equation*}
and
\begin{equation*}
Q(t)=2e^{2t}-e^t \bigl(t^2-6 t+18\bigr)+8(t+2)
\end{equation*}
on $(0,\infty)$.
\par
By expanding the function $Q(t)$ into power series at $t=0$, we have
$$
Q(t)=\sum_{k=3}^\infty\frac{2^{k+1}-(k-9)(k+2)}{k!}t^k>0,\quad t>0.
$$
Direct differentiation gives
\begin{align*}
\frac{\td}{\td t}\biggl[\frac{P(t)}{Q(t)}\biggr]=-\frac{\theta(t)} {(e^t-1)^2[Q(t)]^2},
\end{align*}
where
\begin{align*}
\theta(t)&=2 e^{5 t} \bigl(t^3-15 t^2+78 t-140\bigr)\\
&\quad+e^{4 t} \bigl(t^4-16 t^3+132 t^2-420 t+1236\bigr)\\
&\quad-2 e^{3 t} \bigl(5 t^4-26 t^3+125 t^2-132 t+1018\bigr)\\
&\quad+e^{2 t} \bigl(17 t^4-88 t^3+356 t^2-228 t+1660\bigr)\\
&\quad-e^t \bigl(8 t^4-38 t^3+40 t^2-228 t+756\bigr)+176
\end{align*}
for $t>0$. Straightforward differentiating leads to
\begin{align*}
\theta'(t)&=2 e^t \bigl[e^{4 t} \bigl(5 t^3-72 t^2+360 t-622\bigr)\\
&\quad+2 e^{3 t} \bigl(1131-354 t+120 t^2-15 t^3+t^4\bigr)\\
&\quad-e^{2 t} \bigl(2922-146 t+297 t^2-58 t^3+15 t^4\bigr)\\
&\quad+e^t \bigl(1546+128 t+224 t^2-54 t^3+17 t^4\bigr)\\
&\quad-264+74 t+37 t^2+3 t^3-4 t^4\bigr]\\
&\triangleq2 e^t\theta_1(t),\\
\theta_1'(t)&=e^{4 t} \bigl(20 t^3-273 t^2+1296 t-2128\bigr)\\
&\quad+2 e^{3 t} \bigl(3039-822 t+315 t^2-41 t^3+3 t^4\bigr)\\
&\quad-2 e^{2 t} \bigl(2849+151 t+210 t^2-28 t^3+15 t^4\bigr)\\
&\quad+e^t \bigl(1674+576 t+62 t^2+14 t^3+17 t^4\bigr)\\
&\quad+74+74 t+9 t^2-16 t^3,\\
\theta_1''(t)&=2 e^{4 t} \bigl(40 t^3-516 t^2+2319 t-3608\bigr)\\
&\quad+6 e^{3 t} \bigl(2765-612 t+274 t^2-37 t^3+3 t^4\bigr)\\
&\quad-2 e^{2 t} \bigl(5849+722 t+336 t^2+4 t^3+30 t^4\bigr)\\
&\quad+e^t \bigl(2250+700 t+104 t^2+82 t^3+17 t^4\bigr)\\
&\quad+74+18 t-48 t^2,\\
\theta_1^{(3)}(t)&=2 e^{4 t} \bigl(160 t^3-1944 t^2+8244 t-12113\bigr)\\
&\quad+6 e^{3 t} \bigl(7683-1288 t+711 t^2-99 t^3+9 t^4\bigr)\\
&\quad-8 e^{2 t} \bigl(3105+529 t+171 t^2+32 t^3+15 t^4\bigr)\\
&\quad+e^t \bigl(2950+908 t+350 t^2+150 t^3+17 t^4\bigr)\\
&\quad +18-96 t,\\
\theta_1^{(4)}(t)&=32 e^{4 t} \bigl(40 t^3-456 t^2+1818 t-2513\bigr)\\
&\quad+6 e^{3 t} \bigl(21761-2442 t+1836 t^2-261 t^3+27 t^4\bigr)\\
&\quad-8 e^{2 t} \bigl(6739+1400 t+438 t^2+124 t^3+30 t^4\bigr)\\
&\quad+e^t \bigl(3858+1608 t+800 t^2+218 t^3+17 t^4\bigr)-96,\\
\theta_1^{(5)}(t)&=e^t \bigl[64 e^{3 t} \bigl(80 t^3-852 t^2+3180 t-4117\bigr)\\
&\quad+18 e^{2 t} \bigl(20947-1218 t+1575 t^2-225 t^3+27 t^4\bigr)\\
&\quad-16 e^t \bigl(7439+1838 t+624 t^2+184 t^3+30 t^4\bigr)\\
&\quad+5466+3208 t+1454 t^2+286 t^3+17 t^4\bigr]\\
&\triangleq e^t \theta_2(t),\\
\theta_2'(t)&=2 \bigl[96 e^{3 t} \bigl(80 t^3-772 t^2+2612 t-3057\bigr)\\
&\quad+9 e^{2 t} \bigl(40676+714 t+2475 t^2-342 t^3+54 t^4\bigr)\\
&\quad-8 e^t \bigl(9277+3086 t+1176 t^2+304 t^3+30 t^4\bigr)\\
&\quad+1604+1454 t+429 t^2+34 t^3\bigr],\\
\theta_2''(t)&=4 \bigl[48 e^{3 t} \bigl(240 t^3-2076 t^2+6292 t-6559\bigr)\\
&\quad+9 e^{2 t} \bigl(41033+3189 t+1962 t^2-234 t^3+54 t^4\bigr)\\
&\quad-4 e^t \bigl(12363+5438 t+2088 t^2+424 t^3+30 t^4\bigr)\\
&\quad+727+429 t+51 t^2\bigr],\\
\theta_2^{(3)}(t)&=4 \bigl[48 e^{3 t} \bigl(720 t^3-5508 t^2+14724 t-13385\bigr)\\
&\quad+9 e^{2 t} \bigl(85255+10302 t+3222 t^2-252 t^3+108 t^4\bigr)\\
&\quad-4 e^t \bigl(17801+9614 t+3360 t^2+544 t^3+30 t^4\bigr)\\
&\quad+429+102 t\bigr],\\
\theta_2^{(4)}(t)&=8 \bigl[72 e^{3 t} \bigl(720 t^3-4788 t^2+11052 t-8477\bigr)\\
&\quad+18 e^{2 t} \bigl(45203+6762 t+1422 t^2-18 t^3+54 t^4\bigr)\\
&\quad-2 e^t \bigl(27415+16334 t+4992 t^2+664 t^3+30 t^4\bigr)+51\bigr],\\
\theta_2^{(5)}(t)&=16 e^t \bigl[108 e^{2 t} \bigl(720 t^3-4068 t^2+7860 t-4793\bigr)\\
&\quad+18 e^t \bigl(48584+8184 t+1395 t^2+90 t^3+54 t^4\bigr)\\
&\quad-43749-26318 t-6984 t^2-784 t^3-30 t^4\bigr]\\
&\triangleq16 e^t \theta_3(t),\\
\theta_3'(t)&=2 \bigl[108 e^{2 t} \bigl(720 t^3-2988 t^2+3792 t-863\bigr)\\
&\quad+9 e^t \bigl(56768+10974 t+1665 t^2+306 t^3+54 t^4\bigr)\\
&\quad-13159-6984 t-1176 t^2-60 t^3\bigr],\\
\theta_3''(t)&=6 \bigl[72 e^{2 t} \bigl(1033+804 t-1908 t^2+720 t^3\bigr)\\
&\quad+3 e^t \bigl(67742+14304 t+2583 t^2+522 t^3+54 t^4\bigr)\\
&\quad-4 \bigl(582+196 t+15 t^2\bigr)\bigr],\\
\theta_3^{(3)}(t)&=6 \bigl[144 e^{2 t} \bigl(1435-1104 t-828 t^2+720 t^3\bigr)\\
&\quad+3 e^t \bigl(82046+19470 t+4149 t^2+738 t^3+54 t^4\bigr)\\
&\quad-8 (98+15 t)\bigr],\\
\theta_3^{(4)}(t)&=18 \bigl[96 e^{2 t} \bigl(720 t^3+252 t^2-1932 t+883\bigr)\\
&\quad+e^t \bigl(54
   t^4+954 t^3+6363 t^2+27768 t+101516\bigr)-40\bigr]\\
&>18 \bigl[96 e^{2 t} \bigl(720 t^3+252 t^2-1932 t+883\bigr)\\
&\quad+(1+t)\bigl(54
   t^4+954 t^3+6363 t^2+27768 t+101516\bigr)-40\bigr]\\
&=e^{2 t} \bigl[e^{-2 t} \bigl(101476+129284 t+34131 t^2+7317 t^3+1008 t^4+54 t^5\bigr)\\
&\quad+96 \bigl(883-1932 t+252 t^2+720 t^3\bigr)\bigr].
\end{align*}
\par
In the light of Descartes' Sign Rule, the function
$$
u(t)=883-1932 t+252 t^2+720 t^3
$$
has at most two zeros on $[0,\infty)$. Since $u(0)=883$, $u(1)=-77$ and $u(2)=3787$, these two zeros are all less than $2$, which implies that the function $u(t)$ is positive on $[2,\infty)$, and so $\theta_3^{(4)}(t)>0$ on $[2,\infty)$.
\par
In~\cite[p.~269, 3.6.6]{mit} and~\cite{exp-beograd}, it was listed that
\begin{equation}
e^x\le\frac{2+x}{2-x},\quad 0\le x<2.
\end{equation}
Hence, for $t\in(0,2)$, we have
\begin{align*}
\theta_3^{(4)}(t)&\ge e^{2 t} \biggl[\biggl(\frac{2-t}{2+t}\biggr)^2\bigl(101476+129284 t+34131 t^2+7317 t^3\\
&\quad+1008 t^4+54 t^5\bigr)+96 \bigl(883-1932 t+252 t^2+720 t^3\bigr)\biggr]\\
&=\frac{e^{2 t}}{(t+2)^2}\bigl(54 t^7+792 t^6+72621 t^5+309567 t^4\\
&\quad+209804 t^3-839488 t^2-291584t+744976\bigr)\\
&=\frac{e^{2 t}}{(t+2)^2}\bigl[54 t^7+792 t^6+72621 t^5\\
&\quad +\bigl(309567 t^2+828938 t+508821\bigr) (t-1)^2-102880 t+236155\bigr]\\
&>\frac{236155-102880 t}{(t+2)^2}e^{2 t}\\
&>0.
\end{align*}
In conclusion, the function $\theta_3^{(4)}(t)$ is positive on $(0,\infty)$.
\par
A direct calculation yields
\begin{align*}
  \theta(0)&=0, & \theta'(0)&=0, & \theta_1'(0)&=0, \\
  \theta_1''(0)&=0, & \theta_1^{(3)}(0)&=0, & \theta_1^{(4)}(0)&=0, \\
  \theta_1^{(5)}(0)&=0, & \theta_2(0)&=0, &  \theta_2'(0)&=0, \\
  \theta_2''(0)&=22960, & \theta_2^{(3)}(0)&=216160, & \theta_2^{(4)}(0)&=1188248, \\
  \theta_2^{(5)}(0)&=5009904, & \theta_3(0)&=313119, & \theta_3'(0)&=809098, \\
  \theta_3''(0)&=1651644, &  \theta_3^{(3)}(0)&=2711964, & \theta_3^{(4)}(0)&=3352392.
\end{align*}
This implies that
\begin{align*}
  \theta(t)&>0, & \theta'(t)&>0, & \theta_1'(t)&>0, \\
  \theta_1''(t)&>0, & \theta_1^{(3)}(t)&>0, & \theta_1^{(4)}(t)&>0, \\
  \theta_1^{(5)}(t)&>0, & \theta_2(t)&>0, &  \theta_2'(t)&>0, \\
  \theta_2''(t)&>22960, & \theta_2^{(3)}(t)&>216160, & \theta_2^{(4)}(t)&>1188248, \\
  \theta_2^{(5)}(t)&>5009904, & \theta_3(t)&>313119, & \theta_3'(t)&>809098, \\
  \theta_3''(t)&>1651644, &  \theta_3^{(3)}(t)&>2711964, & \theta_3^{(4)}(t)&>3352392
\end{align*}
on $(0,\infty)$. As a result, the function $\frac{P(t)}{Q(t)}$ is decreasing on $(0,\infty)$, with
$$
\lim_{t\to0}\frac{P(t)}{Q(t)}=1\quad \text{and}\quad \lim_{t\to\infty}\frac{P(t)}{Q(t)}=0.
$$
Hence,
\begin{enumerate}
\item
when $\lambda\le0$, the function $h_\lambda(x)$ is completely monotonic on $(0,\infty)$;
\item
when $\lambda\ge4$, the negative of $h_\lambda(x)$ is completely monotonic on $(0,\infty)$.
\end{enumerate}
Since $\frac2{x^2}$ is completely monotonic on $(0,\infty)$ and the product of finitely many completely monotonic functions is also completely monotonic,
\begin{enumerate}
\item
when $\lambda\le0$, the difference $f_\lambda(x)-f_\lambda(x+1)$ is completely monotonic on $(0,\infty)$, that is,
$$
0\le(-1)^i[f_\lambda(x)-f_\lambda(x+1)]^{(i)}=(-1)^i[f_\lambda(x)]^{(i)}-(-1)^i[f_\lambda(x+1)]^{(i)},\quad i\ge0.
$$
By virtue of induction, we obtain
\begin{multline*}
(-1)^i[f_\lambda(x)]^{(i)}\ge(-1)^i[f_\lambda(x+1)]^{(i)}\ge(-1)^i[f_\lambda(x+2)]^{(i)}\ge\dotsm\\*
\ge(-1)^i[f_\lambda(x+k)]^{(i)} \ge\lim_{k\to\infty}(-1)^i[f_\lambda(x+k)]^{(i)}=0
\end{multline*}
for $i\ge0$. So the function $f_\lambda(x)$ for $\lambda\le0$ is completely monotonic on $(0,\infty)$.
\item
when $\lambda\ge4$, a similar argument leads to the complete monotonicity of $f_\lambda(x)$ on $(0,\infty)$.
\end{enumerate}
The sufficiency is proved.
\par
Multiplying by $x^n$ on both sides of \eqref{psisymp4} yields
\begin{equation}\label{lim-x-2-0}
\lim_{x\to0^+}\bigl[{x^n}\psi^{(n-1)}(x)\bigr]={(-1)^{n}(n-1)!},\quad n\in\mathbb{N}.
\end{equation}
Using L'H\^ospital's rule,  the limit~\eqref{lim-x-2-0}, and the formula~\eqref{psisymp4}, we have
\begin{align*}
&\quad\lim_{x\to0^+}\frac{12x^4(x+1)^2\bigl\{[\psi'(x)]^2+\psi''(x)\bigr\}-x^2-12}x\\
&=\lim_{x\to0^+}\bigl\langle24 (x+1) x^4 \bigl\{[\psi'(x)]^2+\psi''(x)\bigr\}\\
&\quad+12 (x+1)^2 x^4 \bigl[2\psi'(x) \psi''(x)+\psi^{(3)}(x)\bigr]\\
&\quad+48 (x+1)^2 x^3 \bigl\{[\psi'(x)]^2+\psi''(x)\bigr\}-2 x\bigr\rangle\\
&=24\lim_{x\to0^+}\bigl[x^2\psi'(x)\bigr]^2 +12\lim_{x\to0^+}\bigl[x^4\psi^{(3)}(x)\bigr]+48\lim_{x\to0^+}\bigl[x^3\psi''(x)\bigr]\\
&\quad+24\lim_{x\to0^+}\bigl[x^4\psi''(x)\bigr]+24\lim_{x\to0^+}\bigl\{x^4\psi'(x) \psi''(x)+2x^3 [\psi'(x)]^2\bigr\}\\
&=24\lim_{x\to0^+}\bigl\{x^2\psi'(x)\bigl[x^2\psi''(x)+2x\psi'(x)\bigr]\bigr\}\\
&=24\lim_{x\to0^+}\bigl[x^2\psi''(x)+2x\psi'(x)\bigr]\\
&=24\lim_{x\to0^+}\biggl\{x^2\biggl[\psi''(x)+\frac{(-1)^{3-1}(3-1)!}{x^3}\biggr] \\
&\quad+2x\biggl[\psi'(x)+\frac{(-1)^{2-1}(2-1)!}{x^2}\biggr]\biggr\}\\
&=24\lim_{x\to0^+}\bigl[x^2\psi''(x+1)+2x\psi'(x+1)\bigr]\\
&=0.
\end{align*}
\par
In \cite[p.~260, 6.4.12 and 6.4.13]{abram}, it was listed that
\begin{equation}\label{psi'(z).sim.infty}
\psi'(z)\sim\frac1z+\frac1{2z^2}+\frac1{6z^3}-\frac1{30z^5}+\frac1{42z^7}-\frac1{30z^9}+\dotsm
\end{equation}
and
\begin{equation}\label{psi''(z).sim.infty}
\psi''(z)\sim-\frac1{z^2}-\frac1{z^3}-\frac1{2z^4}+\frac1{6z^6}-\frac1{6z^8}+\frac3{10z^{10}} -\frac5{6z^{12}}+\dotsm
\end{equation}
as $z\to\infty$ in $|\arg z|<\pi$. Therefore, we have
\begin{equation*}
\frac{12x^4(x+1)^2\bigl\{[\psi'(x)]^2+\psi''(x)\bigr\}-x^2-12}x\sim4-\frac{82}{15 x}+\frac{14}{3 x^2}-\frac{29}{35 x^3}+\dotsm
\end{equation*}
as $x\to\infty$. Hence
\begin{equation*}
\lim_{x\to\infty}\frac{12x^4(x+1)^2\bigl\{[\psi'(x)]^2+\psi''(x)\bigr\}-x^2-12}x=4.
\end{equation*}
If the function $f_\lambda(x)$ is completely monotonic on $(0,\infty)$, then
$$
[\psi'(x)]^2+\psi''(x)-\frac{x^2+\lambda x+12}{12x^4(x+1)^2}>0
$$
which can be rearranged as
\begin{equation*}
\lambda<\frac{12x^4(x+1)^2\bigl\{[\psi'(x)]^2+\psi''(x)\bigr\}-x^2-12}x\to0
\end{equation*}
as $x\to0^+$, which means that $\lambda\le0$. If $-f_\lambda(x)$ is completely monotonic on $(0,\infty)$, then
\begin{equation*}
\lambda>\frac{12x^4(x+1)^2\bigl\{[\psi'(x)]^2+\psi''(x)\bigr\}-x^2-12}x\to4
\end{equation*}
as $x\to\infty$, so $\lambda\ge4$. The necessity is proved.
\par
The double inequality~\eqref{x-4-upper-improve-ineq} and its best possibility follow from the necessary and sufficient conditions for the function $f_\lambda(x)$ to be completely monotonic on $(0,\infty)$. The proof of Theorem~\ref{x-4-upper-improve-thm} is complete.
\end{proof}

\begin{proof}[Second proof for the first part of Theorem~\ref{x-4-upper-improve-thm}]
It is easy to see that
\begin{equation}
f_\lambda(x)=f(x)-\lambda h(x),
\end{equation}
where $f(x)$ is defined by~\eqref{x-4-upper-f(x)} and
\begin{equation}
  h(x)=\frac1{12x^3(x+1)^2}.
\end{equation}
From~\cite{x-4-di-tri-gamma-upper-lower-combined.tex} it is known that the function $f(x)$ is completely monotonic on $(0,\infty)$, as mentioned on page~\pageref{x-4-upper-f(x)}. Since $h(x)$ is completely monotonic on $(0,\infty)$, it follows that the function $f_\lambda(x)$ is completely monotonic on $(0,\infty)$ when $\lambda\le0$.
Utilizing~\eqref{psisymp4} for $n=1$ and $n=2$ yields that
\begin{align*}
f_\lambda(x)&=\biggr[\psi'(x+1)+\frac1{x^2}\biggl]^2 +\biggl[\psi''(x+1)-\frac{2}{x^3}\biggr] -\frac{x^2+\lambda x+12}{12x^4(x+1)^2}\\
&=-\frac{\lambda}{12}\biggl[\frac{1}{x^3}-\frac{2}{x^2}+\frac{3}{x}-\frac{4+3 x}{(1+x)^2}\biggr] -\frac{37}{12 x^2}+\frac{25}{6 x}-\frac{63+50 x}{12 (1+x)^2} \\
&\quad+\frac{2}{x^2}\psi'(x+1)+\bigl[\psi'(x+1)\bigr]^2+\psi''(x+1).
\end{align*}
Therefore, if $\lambda>0$ then $\lim_{x\to0^+}f_\lambda(x)=-\infty$. This implies that the function $f_\lambda(x)$ is completely monotonic on $(0,\infty)$ if and only if $\lambda\le0$.
\end{proof}

\begin{proof}[Second proof for the necessity of the second part of Theorem~\ref{x-4-upper-improve-thm}]
From the asymptotic expansions in~\eqref{psi'(z).sim.infty} and~\eqref{psi''(z).sim.infty} it is not difficult to obtain that
\begin{equation*}
f_\lambda(x)=\frac{4-\lambda}{12x^5}+O\biggl(\frac1{x^6}\biggr)
\end{equation*}
as $x\to\infty$. This implies that the function $f_\lambda(x)$ is positive for small positive number $x$ when $\lambda<4$. So the necessary condition in the second part of Theorem~\ref{x-4-upper-improve-thm} follows.
\end{proof}

\begin{proof}[Second proof for a part of the sufficiency in the second part of Theorem~\ref{x-4-upper-improve-thm}]
It is clear that
$$
-f_\lambda(x)=g(x)+\frac{\lambda-13}{12 x^3 (1+x)^2},
$$
where $g(x)$ is defined by~\eqref{x-4-upper-g(x)}.
From the complete monotonicity of the function $g(x)$ mentioned on page~\pageref{x-4-upper-g(x)}, it follows that if $\lambda\ge13$ the function $-f_\lambda(x)$ is completely monotonic on $(0,\infty)$. This gives an alternative proof for a part of the sufficiency in the second part of Theorem~\ref{x-4-upper-improve-thm}.
\end{proof}

\begin{rem}
The first proof of Theorem~\ref{x-4-upper-improve-thm} is direct and independent, but other partial proofs are based on the main result in~\cite{x-4-di-tri-gamma-upper-lower-combined.tex}.
\end{rem}

\section{More remarks}

In this section we derive some corollaries from the first proof of Theorem~\ref{x-4-upper-improve-thm} and pose a double inequality for bounding the function $[\psi'(x)]^2+\psi''(x)$ on $(0,\infty)$.

\begin{rem}
It is easy to see that the double inequality \eqref{x-4-upper-improve-ineq} for $\mu=0$ and $\nu=4$ recovers the left-hand side inequality and refines the right-hand side inequalities in~\eqref{psi'psi''-upper-weaker} and~\eqref{psi'psi''-upper}.
\par
It is clear that the bounds in \eqref{psi'psi''-upper} and~\eqref{x-4-upper-improve-ineq} are simpler than~\eqref{psi'psi''-lower} and~\eqref{psi'psi''-upper-weaker}.
\end{rem}

\begin{rem}
From the proof of Theorem~\ref{x-4-upper-improve-thm}, it is easy to deduce that the function
\begin{multline}\label{psi'-ineq-lower}
\psi'(x)-\frac{\lambda \bigl(4+8 x+5 x^2\bigr)}{24 x (1+x)^3 (2+x)^2}\\*
-\frac{24+120 x+283 x^2+399 x^3+345 x^4+181 x^5+51 x^6+6 x^7}{6 x^2 (1+x)^4 (2+x)^2}
\end{multline}
is completely monotonic on $(0,\infty)$ if and only if $\lambda\le0$, and so is its negative if and only if $\lambda\ge4$.
\end{rem}

\begin{rem}
Integrating the function~\eqref{psi'-ineq-lower} and considering its positivity, it is immediate to see that the double inequality
\begin{multline}\label{psi-ineq-double-xi-eta}
\xi\biggl[\frac{\ln x-9 \ln(1+x)+8 \ln(2+x)}{24}-\frac{18+33 x+14 x^2}{48 (1+x)^2 (2+x)}\biggr]\\* <\psi(x)-\frac{28 x^4+87 x^3+73 x^2+3 x-12}{6 x (x+1)^3 (x+2)}-\frac{43\ln (x+1)-37\ln(x+2)}{6}\\
<\eta\biggl[\frac{\ln x-9 \ln(1+x)+8 \ln(2+x)}{24}-\frac{18+33 x+14 x^2}{48 (1+x)^2 (2+x)}\biggr]
\end{multline}
holds on $(0,\infty)$ if and only if $\xi\ge4$ and $\eta\le0$.
\end{rem}

\begin{rem}
We note that the function
\begin{equation}
G_\lambda(x)=\ln\Gamma(x)+x+\frac{1}{12 (1+x)^2}+\frac{\lambda-48}{48 (1+x)}+\frac{H_\lambda(x)}{24}
\end{equation}
\begin{enumerate}
\item
is completely monotonic on $(0,\infty)$ if and only if $\lambda\le0$;
\item
satisfies $(-1)^{k+1}[G_\lambda(x)]^{(k)}>0$ for $k\in\mathbb{N}$ on $(0,\infty)$ if and only if $\lambda\ge4$;
\item
has limits
\begin{equation}
\lim_{x\to0^+}G_\lambda(x)=\frac{488 \ln2-44+\lambda(1-24\ln2)}{48}
\end{equation}
and
\begin{equation}
\lim_{x\to\infty}G_\lambda(x)=\frac{124+12 \ln(2\pi)-7\lambda}{24},
\end{equation}
\end{enumerate}
where
\begin{equation}
  \begin{split}
   H_\lambda(x)&=(24-\lambda x) \ln x+[(9\lambda-172)x+12\lambda-256]\ln(x+1)\\
    &\quad+4[(37-2\lambda)x-3\lambda+61]\ln(x+2).
  \end{split}
\end{equation}
In particular, we have
\begin{equation}
\begin{split}
G_0(x)&=\ln\Gamma(x)+\ln x+x+\frac{1}{12 (1+x)^2}-\frac{1}{1+x}\\
&\quad+\frac{(61+37 x) \ln(2+x)-(64+43 x) \ln(1+x)}6
\end{split}
\end{equation}
and
\begin{equation}
  \begin{split}
    G_4(x)&=\ln\Gamma(x)+\ln x+x+\frac{1}{12 (1+x)^2}-\frac{11}{12 (1+x)}\\
    &\quad+\frac{(49+29 x) \ln(2+x)-2 (26+17 x) \ln(1+x)-x \ln x}6
  \end{split}
\end{equation}
on $(0,\infty)$.
\end{rem}

\begin{rem}
We conjecture that the double inequality
\begin{equation}\label{x-4-upper-improve-ineq-conj}
\frac1{x^4}\biggl[\frac{x^2+4x+12}{12(x+1)^2}\biggr]^\alpha<[\psi'(x)]^2+\psi''(x) <\frac1{x^4}\biggl[\frac{x^2+4x+12}{12(x+1)^2}\biggr]^\beta
\end{equation}
holds on $(0,\infty)$ if and only if $\alpha\ge\frac65$ and $\beta\le1$.
\end{rem}

\end{document}